\documentclass[12pt]{article}
\usepackage{graphicx,setspace}
\usepackage{amsmath,amssymb,amsthm}
\newcommand{\be}{\begin{equation}}
\newcommand{\ee}{\end{equation}}
\newcommand{\bea}{\begin{eqnarray}}
\newcommand{\beas}{\begin{eqnarray*}}
\newcommand{\no}{\nonumber}
\newcommand{\eea}{\end{eqnarray}}
\newcommand{\eeas}{\end{eqnarray*}}

\newcommand{\remove}[1]{}

\newtheorem{thm}{Theorem}[section]

\newtheorem{lem}[thm]{Lemma}
\newtheorem{prop}[thm]{Proposition}
\newtheorem{remark}[thm]{\bf Remark}

\def\sec{\setcounter{equation}{0}}

\newcounter{cnt1}

\newcounter{cnt3}
\newcommand{\blr}{\begin{list}{$($\roman{cnt1}$)$}
 {\usecounter{cnt1} \setlength{\topsep}{0pt}
 \setlength{\itemsep}{0pt}}}
\newcommand{\bla}{\begin{list}{$($\betaph{cnt2}$)$}
 {\usecounter{cnt2} \setlength{\topsep}{0pt}
 \setlength{\itemsep}{0pt}}}
\newcommand{\bln}{\begin{list}{$($\arabic{cnt3}$)$}
 {\usecounter{cnt3} \setlength{\topsep}{0pt}
 \setlength{\itemsep}{0pt}}}
\newcommand{\el}{\end{list}}

\newcommand{\ep}{\epsilon}
\newcommand{\al}{\alpha}

\newcommand{\del}{\delta}

\def\ka{\kappa}
\def\Del{\Delta}

\def\th{\theta}
\def\ga{\gamma}

\def\mR{\mathbb{R}}
\def\mZ{\mathbb{Z}}

\def\cF{{\mathcal F}}

\def\Pn{{\mathcal P}_{n}}

\def\P{\mathcal P}

\def\1{\mathbf{1}}
\def\half{\frac{1}{2}}
\def\hr{\hat{r}}
\def\tr{\tilde{r}}
\def\br{\bar{r}}
\def\tL{\tilde{L}}
\def\tE{\tilde{E}}
\def\lb{\left(}
\def\rb{\right)}
\def\f{\frac}
\def\lab{\label}
\def\nn{\nonumber}
\begin{document}

\begin{titlepage}
\begin{center}
{\bf Connecting the Random Connection Model} \\
\vspace{0.2in} {Srikanth K. Iyer \footnote{corresponding author: srikiyer@gmail.com}$^,$\footnote{Research Supported in part by UGC Center for Advanced Studies}}\\
Department of Mathematics\\
Indian Institute of Science, Bangalore, India.
\end{center}
\vspace{0.1in}
%
\sloppy
\begin{center} {\bf Abstract} \end{center}


{Consider the random graph $G(\Pn,r)$ whose vertex set  $\Pn$ is a Poisson point process of intensity $n$ on $(- \half, \half]^d$, $d \geq 2$. Any two vertices $X_i,X_j \in \Pn$ are connected by an edge with probability $g\lb \f{d(X_i,X_j)}{r} \rb$, independently of all other edges, and independent of the other points of $\Pn$. $d$ is the toroidal metric, $r > 0$ and $g:[0,\infty) \to [0,1]$ is non-increasing and $\al = \int_{\mR^d} g(|x|) dx < \infty$. Under suitable conditions on $g$, almost surely, the critical parameter $d_n$ for which $G(\Pn, \cdot)$ does not have any isolated nodes satisfies $\lim_{n \to \infty} \f{\al n d_n^d}{\log n} = 1$. Let $\beta = \inf\{x > 0: x g\lb \f{\al}{x \th} \rb > 1 \}$, and $\th$ be the volume of the unit ball in $\mR^d$. Then for all $\ga > \beta$, $G\lb\Pn, \lb \f{\ga \log n}{\al n} \rb^{\f{1}{d}}\rb$ is connected with probability approaching one as $n \to \infty$. The bound can be seen to be tight for the usual random geometric graph obtained by setting $g = 1_{[0,1]}$. We also prove some useful results on the asymptotic behavior of the length of the edges and the degree distribution in the {\it connectivity regime}.
}\\
%

\vspace{0.1in}
{\sl AMS 1991 subject classifications}: \\
\hspace*{0.5in} Primary:   60D05, 60G70;
\hspace*{0.5in} Secondary:  05C05, 90C27 

\end{titlepage}

\section{Introduction and Main Results}
\label{sec:intro}
\sec
The random connection model (RCM) has proved to be a very useful model in applications in many branches of science such as physics, epidemiology and telecommunications \cite{Meester07}. A specific case of the RCM is the random geometric graph (RGG) studied in great detail in \cite{Penrose03}. The RGG, the RCM, and its generalizations such as the signal-to-interference-plus-noise ratio (SINR) graph have been used extensively in modeling in wireless networks (\cite{Meester07, Haenggi13}). In this context, some of the questions of interest pertain to percolation, connectivity and coverage. The existence of a phase transition for the percolation problem in a RCM was shown in \cite{Penrose91}. Random connection models have also been studied in \cite{Ballister04, Ballister05}. In recent years there is a growing interest in the study of spatial preferential attachment graphs. Such a graph incorporates elements of the RCM as well as those of preferential attachment graphs. The problem of interest is  the asymptotic degree distribution, clustering coefficient and phase transition behavior (see for example \cite{Jacob15} and references therein). However, the problem of connectivity in a RCM still remains open. \cite{Mao11} derive a parameter regime under which the number of isolated nodes in the RCM converges to a Poisson distribution. This provides a necessary condition for the RCM to be connected with probability approaching one (also referred to as {\it with high probability ( whp)}). Recently \cite{Penrose15} considers a general random connection model and derives conditions under which one obtains Poisson convergence for the number of isolated vertices. In this paper we derive a strong law result for the critical parameter required to eleminate isolated nodes in a RCM. This critical parameter can be thought of as the analog of the largest nearest neighbor distance in the usual random geometric graph (see \cite{Penrose03}). The main result of this paper is a sufficient condition for the RCM to be connected {\it whp}. We also derive strong law results for the length of the longest edge in the connectivity regime. These results are of independent interest, for instance, in determining the diameter of the graph and in routing algorithms in communication networks. In this paper they enable us to work a localization argument in a couple of results. 

In order to describe some of our results it will be useful to construct a coupled family of graphs. Let $(\Omega, \cF, P)$ be a probability space on which all the random variables to be described are defined. Let $S := \lb -\half, \half \right]^d$, where $d \geq 2$. We ignore the edge effects in our graphs by equipping $S$ with the toroidal metric $d(x,y) = \inf\{|x-y-z|: z \in \mZ^d\}$, where $| \cdot |$ denotes the Euclidean norm. Let $\{X_1, X_2, \ldots \}$ be a sequence of independent random variables distributed uniformly in $S$. Let $\{N_n\}_{n \geq 1}$ be a non-decreasing sequence of random variables with $N_n$ having a Poisson distribution with mean $n$. Define the coupled sequence of Poisson point processes on $S$ by $\Pn := \{X_1, X_2, \ldots , X_{N_n}\},$ $n \geq 1$. The connection function $g:[0,\infty) \to [0,1]$ is a non-increasing function satisfying $\al = \int_{\mR^d} g(|x|) dx < \infty$. Let $\{U_{ij}\}_{i,j \geq 1}$ be a sequence of independent random variables with uniform distribution on $[0,1]$. The graph $G_n(r) = G(\Pn,r)$ is the graph with vertex set $\Pn$ with an edge between two vertices $X_i,X_j \in \Pn$ provided $U_{ij} \leq g\lb \f{d(X_i,X_j)}{r} \rb$, independent of every thing else.
Thus an edge exists between $X_i$ and $X_j$ in $G_n(r)$ with probability $g\lb \f{d(X_i,X_j)}{r} \rb$ which is non-decreasing in $r$. This together with the coupling constructed above ensures that if $X_i,X_j \in \Pn$ and an edge exists between these two vertices in $G_{n+1}(r)$ then it exists between these two vertices in $G_n(s)$ for all $s > r$. With $g = 1_{[0,1]}$, the RCM reduces to the usual RGG. Throughout this paper, a coupling will refer to the construction above.

For RGG it is known that the critical radius required to connect the graph and that of the largest nearest neighbor distance (the smallest radius required to eliminate isolated nodes) are identical {\it whp}. As one increases the radius in a RGG, with high probability,  it gets connected at the same time when the isolated nodes disappear (see (13.37) \cite{Penrose03}). Thus to study connectivity in a RCM it is instructive to consider the regime under which the expected number of isolated nodes stabilize. Let $W_n(r)$ be the number of isolated nodes (vertices with degree zero) in $G_n(r)$.

\begin{prop} Define the sequence of parameters $\{r_n(b)\}_{n \geq 1}$ by $r_n(b)^d = \f{\log n + b}{\al n}$, $b \in \mR$ and suppose that the connection function satisfies $g(r) = o(r^{-c})$ as $r \to \infty$ for some $c > d$. Then $E[W_n(r_n(b))] \to e^{-b}$ as $n \to \infty$.
\lab{p1}
\end{prop}

The result also holds if we replace $b$ by a sequence $b_n \to b$. Mao and Anderson (2011) in fact show that with $r_n(b)$ as in Proposition \ref{p1}, the distribution of $W_n$ converges to a Poisson distribution with mean $e^{-b}$. This shows that $b_n \to \infty$ is a necessary condition for $G_n(r_n(b))$ to be connected {\it whp}. To state the strong law result for the critical parameter required to eleminate isolated nodes, we define for $n \geq 1$ the sequence of random variables
\be d_n = \inf\{r > 0: W_n(r) = 0\}. \lab{e1} \ee
\begin{thm} $(i)$ Suppose that $g(r) = o(r^{-c})$ as $r \to \infty$ for some $c > d$. Then almost surely
\be \limsup_{n \to \infty} \f{\al n d_n^d}{\log n} \leq 1. \lab{e2} \ee
$(ii)$ If $g(r) = o(e^{-c r})$ as $r \to \infty$ for some $c > 0$, then almost surely
\be \lim_{n \to \infty} \f{\al n d_n^d}{\log n} = 1. \lab{e3} \ee
$(iii)$ If $g(r) = o(r^{-c})$ as $r \to \infty$ for some $c > 3d$, then almost surely
\be \liminf_{n \to \infty} \f{\al n d_n^d}{\log n} \geq \f{c - 3d}{c - d}. \lab{e3a} \ee
\lab{t1}
\end{thm}

We are now ready to state the main result of this paper which is a sufficient condition for the graph to be connected {\it whp}. The idea behind the proof is similar to the renormalization argument used in percolation problems. Let $\th$ be the volume of the unit ball in $\mR^d$. Since the function $g$ is non-increasing and takes values in the interval $[0,1]$, the function $x g\lb \f{\al}{x \th} \rb$ is strictly increasing with limits at zero and infinity being zero and infinity respectively. Thus 
\be
\beta =\inf\{ x > 0:  x g\lb \f{\al}{x \th} \rb > 1 \},
\lab{e3b}
\ee
is well defined and $\beta \geq 1$. For any $\ga > 0$ and $n \geq 1$, $\hr_n(\ga)$ is defined by the following equation. 
\be \hr(\ga)^d =  \f{\ga \log n}{\al n}.
\lab{e4}
\ee
A sequence of events $\{A_n\}_{n \geq 1}$ is said to hold {\it whp} if $P(A_n) \to 1$ as $n \to \infty$. 

\begin{thm} Suppose that the connection function satisfies $g(r) = o(r^{-c})$ as $r \to \infty$ for some $c > d$. Let $\beta$ be as defined in (\ref{e3b}) and $\hr_n$ as in (\ref{e4}). Then for any $\ga > \beta$ the sequence of graphs $G_n\lb\hr_n(\ga)\rb$ is connected {\it whp}.
\lab{t2}
\end{thm}

\begin{remark} For the usual random geometric graph obtained by setting $g = 1_{[0,1]}$ we see that $\beta = 1$. Thus the bound obtained in Theorem \ref{t2} reduces to the condition required to connect the RGG {\it whp}. The proof of  Theorem \ref{t2} uses the connectivity threshold result for a RGG and thus does not give a new proof of connectivity for RGG.
\end{remark}

\begin{remark} The results of this paper hold true when the number of vertices is fixed to be $n$ instead of the Poisson number $N_n$ by the standard de-Poissoinization argument (see \cite{Penrose03}). 
\end{remark}
%

\section{Edge lengths and Degree Distributions}
\lab{sec:degree}
\sec

In this section we first prove two results, one on the length of longest edge at a typical node and the other on the longest edge in the graph in the connectivity regime, that is when the scaling parameter satisfies 
%
%
(\ref{e4}) for some $\ga > 0$. In (\ref{e4}) and in what follows, all definitions will be assumed to hold for $n$ sufficiently large so that $\hr_n^d > 0$. We also prove strong law results for the minimum and maximum vertex degrees in the connectivity regime. Given a network, the diameter and the degree distribution are important characteristic that are easy to measure. Thus they can be used as a criterion for choosing a model.  Unlike the usual RGG where the edges have length smaller than the cutoff radius, in the random connection model one can, in principle have edges of arbitrary length. The presence of sufficient number of long edges, for instance, can help reduce the number of ``hops'' required for a packet to reach the destination in a multi-hop wireless network. 

We first consider the length of the longest edge incident on a typical edge. Denote the origin by $O$ and let $P^o$ denote the Palm distribution of $\Pn$ conditioned on a point being located at the origin. Since $\Pn$ is a Poisson point process, the distribution of $\Pn$ under $P^o$ is the same as that of $\Pn \cup \{O\}$ under $P$. Let $\hr_n(\ga)$ be as defined in (\ref{e4}). We will write $\hr_n$ for $\hr_n(\ga)$ whenever there is no ambiguity. Let $B(x,r)$ denote a ball of radius $r$ centered at $x$ with respect to the Euclidean norm denoted by $| \cdot |$. Define
\be G(a) = \int_{B(O,a)^c} g(|z|) dz.
\lab{e4a}
\ee
Let $\beta \in \mR$. Since $n \hr_n^d \to \infty$ as $n \to \infty$ and $G(a) \to 0$ as $a \to \infty$, there exists a sequence $\{a_n\}$ such that $n \hr_n^d G(a_n) \to e^{-\beta}$. Define the sequence $\br_n =  a_n \hr_n$. Our first result is on the length of the longest edge incident on a typical node in $G_n(\hr_n(\ga))$.

\begin{prop} Let $g(r) = o(r^{-c})$ as $r \to \infty$ for some $c > d$. Let $\hr_n(\ga), G(a)$ be as defined in (\ref{e4}), (\ref{e4a}) respectively. For any $\beta \in \mR$, let $a_n, \br_n$ be the sequences defined above. Let $L_n^o$ be the length of the longest edge incident on the origin in the graph $G_n(r_n)$ under the Palm measure $P^o$. Then
\be P^o \lb L_n^o \leq \br_n \rb \to e^{-e^{-\beta}}, \qquad \mbox{ as } n \to \infty.
\lab{e5}
\ee
\lab{p2}
\end{prop}

\begin{remark} The longest edge of a typical node in the RCM is longer by a factor of $a_n$ than that in a RGG in the {\it connectivity regime}.
\lab{rem1}
\end{remark}

The next result is on the length of the longest edge in the graph $G_n(\hr_n)$, where $\hr_n(\ga)$ is as defined in (\ref{e4}). While this result is of independent interest, a variation of this result will be used to carry out a localization argument in the proof of parts $(ii), (iii)$ of Theorem~\ref{t1}.

\begin{prop} Let $g(r) = o(e^{-cr})$ as $r \to \infty$ for some $c > 0$, and let $\hr_n(\ga)$ be as in (\ref{e4}). Then almost surely, the length of the longest edge $L_n(\ga)$ in the graph $G_n(\hr_n(\ga))$ satisfies
\be \limsup_{n \to \infty} \f{L_n(\ga)}{\hr_n(\ga) \log n} \leq \f{1}{c}.
\lab{pe3}
\ee
If $g(r) = o(r^{-c})$ as $r \to \infty$ for some $c > 2d$, then almost surely 
\be \limsup_{n \to \infty} \f{\log L_n(\ga)}{\log n} \leq - \f{c - 2d}{d(c-d)}.
\lab{pe3a}
\ee
\lab{p3}
\end{prop}

We now state some results on the vertex degrees in $G_n(r_n)$. 

\begin{prop} Let $\hr_n$ be as in (\ref{e4}), $g(r) = o(r^{-c})$ as $r \to \infty$ for some $c > d$. Define $D_n(k)$ to be the number of vertices in $G_n(\hr_n)$ of degree at least $k$. Then
\be
E \left[ \f{D_n(k_n)}{n} \right] \to 1 - \Phi(d),
\lab{e6}
\ee
where $k_n = \al n \hr_n^d + d \sqrt{\al n \hr_n^d}$ and $\Phi$ is the standard normal distribution function.
\lab{p4}
\end{prop}

The next two results give strong law asymptotics for the maximum and minimum vertex degrees. Let $H : [0, \infty) \to [0, \infty)$ be defined by $H(0) = 1$, and 
\be
H(x) = 1 - x - x \log x, \qquad x > 0.
\lab{e7}
\ee

Note that $H(1) = 0$, $H^{\prime}(x) < 0$ for $x < 1$ and $H^{\prime}(x) > 0$ for $x > 1$. Define $H_+^{-1}:[0,\infty) \to [1,\infty)$ to be the inverse of the restriction of $H$ to $[1,\infty)$ and $H_-^{-1} : [0,1] \to [0,1]$ be the inverse of the restriction of $H$ to $[0,1]$.

\begin{thm} Let $\Del_n = \Del(\hr_n)$ denote the maximum vertex degree in $G_n(\hr_n)$ where $\hr_n$ is as defined in (\ref{e4}). 

$(i)$ Suppose that $g(r) = o(r^{-c})$ as $r \to \infty$ for some $c > d$. Then almost surely,
\be
\limsup_{n \to \infty} \f{\Del_n}{ \log n } \leq \ga H_+^{-1} \lb \f{1}{\ga} \rb. 
\lab{e8}
\ee
$(ii)$ If $g(r) = o(e^{-cr})$ as $r \to \infty$ for some $c > 0$, then
\be
\lim_{n \to \infty} \f{\Del_n}{ \log n } = \ga H_+^{-1} \lb \f{1}{\ga} \rb. 
\lab{e9}
\ee
$(iii)$ If $g(r) = o(r^{-c})$ as $r \to \infty$ for some $c > 3d$, then almost surely
\be
\liminf_{n \to \infty} \f{\Del_n}{ \log n } \geq \ga H_+^{-1} \lb \f{c-3d}{\ga(c-d)} \rb. 
\lab{e10}
\ee
\lab{t3}
\end{thm}

Thus, in the connectivity regime, the expected vertex degree as well as the maximum degree grows logarithmically in $n$. The following result shows that when the scaling parameter is above the critical threshold required to eleminate isolated nodes, then the minimum vertex degree also grows logarithmically in $n$.

\begin{thm} Let $\del_n = \del(\hr_n)$ denote the minimum vertex degree in $G_n(\hr_n)$ where $\hr_n$ is as defined in (\ref{e4}). 

$(i)$ If $g(r) = o(e^{-cr})$ for some $c > 0$ as $r \to \infty$, then almost surely $\del_n \to 0$ as $n \to \infty$ for any $\ga < 1$ .

$(ii)$ If $g(r) = o(r^{-c})$ for some $c > d$ as $r \to \infty$, then for any $\ga > 1$ we have almost surely
\be
\liminf_{n \to \infty} \f{\del_n}{\log n} \geq \ga H_-^{-1} \lb \f{1}{\ga} \rb.
\lab{e11}
\ee
$(iii)$ If $g(r) = o(e^{-cr})$ as $r \to \infty$ for some $c > 0$, then for any $\ga > 1$ we have almost surely
\be
\lim_{n \to \infty} \f{\del_n}{ \log n } = \ga H_-^{-1} \lb \f{1}{\ga} \rb. 
\lab{e12}
\ee
$(iii)$ If $g(r) = o(r^{-c})$ as $r \to \infty$ for some $c > 3d$, then 
\be
\limsup_{n \to \infty} \f{\del_n}{ \log n } \leq \ga H_-^{-1} \lb \f{c-3d}{\ga(c-d)} \rb, 
\lab{e13}
\ee
almost surely for any $ \ga > \f{c-3d}{c-d}$.
\lab{t4}
\end{thm}

\section{Proofs}
\label{sec:intro}
\sec

In what follows $C_1, C_2, \ldots$ will denote constants whose value will change from place to place. We recall a few notations. $O$ denotes the origin in $\mR^d$ and $B(x,r)$ the ball of radius $r$ centered at $x$. $P^o$ denotes the Palm distribution of $\Pn$, that is the measure $P$ conditioned on there being a point of $\Pn$ at the origin and $E^o$ the expectation with respect to $P^o$. We begin the proofs with a useful Lemma that will be invoked several times. 

\begin{lem} Let $\{r_n\}_{n \geq 1}$ be a sequence defined by $r_n^d = \f{a \log n + b}{\al n},$ for some $a > 0$ and $b \in \mR$. If $g(r) = o(r^{-c})$ as $r \to \infty$ for some $c > d$, then 
\be nr_n^d \int_{\mR^d \setminus r_n^{-1}S} g(|z|) dz \to 0, \qquad n \to \infty.
\lab{ep1} \ee
\lab{l1}
\end{lem}

{\bf Proof. } Consider the integral in (\ref{ep1}).
\begin{eqnarray*}
\int_{\mR^d \setminus r_n^{-1}S} g(|z|) dz & \leq & \int_{\mR^d \setminus B(O, r_n^{-1}/2)} g(|z|) dz \\
 & = & C_1 \int_{r_n^{-1}/2}^{\infty} r^{d-1} g(r) dr \\
 & \leq & C_2 \int_{r_n^{-1}/2}^{\infty} r^{d-1-c} dr = C_3 r_n^{c-d}.
\end{eqnarray*}
The result now follows since $c > d$.
\qed

{\bf Proof of Proposition \ref{p1}.} By the Campbell-Mecke formula, $E[W_n(r_n(b)]$ equals $n$ times the probability that the origin is isolated in $G_n(r_n(b))$ under $P^o$. The set of points of $\Pn$ with edges incident on $O$ form a Poisson point process with intensity $n g \lb \f{d(O,\cdot)}{r_n} \rb$ over $S$. Consequently
\[ E[W_n(r_n(b))]  =  n \exp \lb - n \int_S g \lb \f{d(O,y)}{r_n} \rb dy \rb. \]
Note that $d(O,y) = |y|$, $y \in S$. Making the change of variable $z = r_n^{-1} y$ we get
\begin{eqnarray} 
E[W_n(r_n(b))]  & = &  n \exp \lb - n r_n^d \int_{r_n^{-1}S} g(|z|) dz \rb \nn \\
 & = &  n \exp \lb - n r_n^d \lb \al - \int_{\mR^d \setminus r_n^{-1}S} g(|z|) dz \rb \rb \nn \\
 & = & e^{-\beta} \exp \lb n r_n^d \int_{\mR^d \setminus r_n^{-1}S} g(|z|) dz \rb  \to e^{-\beta}.
\lab{pe2}
\end{eqnarray}
as $n \to \infty$ by Lemma \ref{l1}. \qed

The proof of Theorem~\ref{t1} has several parts and so we split it into two parts.

\begin{prop} Let $d_n$ be as defined in (\ref{e1}). 
Suppose that $g(r) = o(r^{-c})$ as $r \to \infty$ for some $c > d$. Then almost surely,
\[ \limsup_{n \to \infty} \f{\al n d_n^d}{\log n} \leq 1. \]
\lab{p5}
\end{prop}

{\bf Proof. } Fix $b>1$ and choose $a$ such that $a(b-1) > 1$. Let $\hr_n$ be as defined in (\ref{e4}). For $k \geq 1$ define $n_k = a^k$ and define the two sequences of events
\[ A_n = \{W_n(\hr_n(b)) > 0 \}, \qquad n \geq 1,\]
and
\[ B_k = \cup_{n=n_k}^{n_{k+1}} A_n, \qquad k \geq 1. \]
The coupling used in the construction of $G_n$ yields
\beas 
P(B_k) & \leq & P\lb W_{n_{k+1}}(\hr_{n_k}(b)) > 0 \rb \\
 & \leq & E \left[ W_{n_{k+1}}(\hr_{n_k}(b)) \right]. 
\eeas
Applying the Campbell-Mecke formula to the right hand side of the above inequality we get
\beas 
P(B_k) & \leq & n_{k+1} P^o \lb O \mbox{ is isolated in } G_{n_{k+1}}(\hr_{n_k}(b)) \rb \\
 & = & n_{k+1} \exp \lb - n_{k+1} \int_S g \lb \f{|y|}{\hr_{n_k}(b)} \rb dy \rb .
\eeas
Make the standard change of variable $z = \hr_{n_k}(b)^{-1}y$. Since $\f{n_{k+1}}{n_k} > 1$ and bounded we have using Lemma~\ref{l1} 
\beas 
P(B_k) & \leq & C_1 n_{k+1} \exp \lb - n_{k+1} \hr_{n_k}(b)^d \al \rb \\
 & \leq & C_2 n_{k+1} \exp \lb - \f{n_{k+1}}{n_k} b \log n_k \rb \\
 & \leq &  \f{C_3}{n_k^{b-1}} = \f{C_3 }{k^{a(b-1)}},
\eeas
which is summable since $a(b-1) > 1$. It follows by the Borel-Cantelli lemma that almost surely only finitely many of the events $B_k$, and hence the events $A_n$, happen. In other words, almost surely
\[ d_n \leq \hr_n(b) \qquad \mbox{or} \qquad \f{\al n d_n^d}{\log n} \leq b, \]
eventually for all $b > 1$. The result now follows since $b > 1$ is arbitrary. \qed

At this stage we will prove Proposition~\ref{p3} since as noted earlier, a variation of this result (see remark below following the proof) will be used to prove parts $(ii), (iii)$ of Theorem~\ref{t1}.

{\bf Proof of Proposition \ref{p3}. } We first prove (\ref{pe3}). Fix $\ep > 0$ and choose $\eta > 0$ sufficiently small so that $(1+\ep)(1-\eta) > 1$. Define the subsequence $n_k = k^a,$ $k \geq 1$, where $a > 0$ is chosen so as to satisfy $a((1+\ep)(1-\eta) - 1) > 1$. Let $D_n(r,L)$ be the event that there is an edge in $G_n(r)$ of length larger than $L$. Define
\[ \tL_n(\ep) = c^{-1} (1 + \ep) \hr_n(\ga) \log n, \]
and let $\{E_k\}_{k \geq 1}$ be the sequence of events defined by
\[ E_k = \cup_{n=n_k}^{n_{k+1}} D_n(\hr_n(\ga), \tL_n(\ep)). \]
By the remark following the construction of the coupled family of graphs $G_n(r)$, we get
\[ E_k \subset D_{n_{k+1}}(\hr_{n_k}(\ga), \tL_{n_{k+1}}(\ep)). \]
Let $A_k$ be the set of points in $\P_{n_{k+1}}$ for which there is an edge of length greater than $\tL_{n_{k+1}}(\ep)$ incident on it in the graph $G\lb \P_{n_{k+1}}, \hr_{n_k}(\ga) \rb$. Using the above inclusion, we get
\[ P(E_k) \leq E \left[ \sum_{X \in \P_{n_{k+1}}} 1_{A_k}(X) \right]. \]
By the Campbell-Mecke formula and using the fact that $1 - e^{-x} \leq x$, we get 
\begin{eqnarray}
P(E_k) & \leq & n_{k+1} P^o(O \in A_k) \nn \\
 & = & n_{k+1} \lb 1 - \exp \lb - n_{k+1} \int_{S \setminus B(O, \tL_{n_{k+1}}(\ep))} g \lb \f{|z|}{\hr_{n_k}(\ga)} \rb dz \rb \rb \nn \\
 & \leq & n_{k+1}^2 \hr_{n_k}(\ga)^d \int_{\hr_{n_k}(\ga)^{-1} S \setminus B(O, \hr_{n_k}(\ga)^{-1} \tL_{n_{k+1}}(\ep))} g \lb |z| \rb dz \nn \\
 & \leq & n_{k+1}^2 \hr_{n_k}(\ga)^d \int_{\mR^d \setminus B(O, \hr_{n_k}(\ga)^{-1} \tL_{n_{k+1}}(\ep))} g \lb |z| \rb dz \nn \\
 & \leq & C_1 n_{k+1}^2 \hr_{n_k}(\ga)^d \int_{\hr_{n_k}(\ga)^{-1} \tL_{n_{k+1}}(\ep)}^{\infty} r^{d-1}g(r) dr.
\lab{pe4}
\end{eqnarray}
By the assumption on the function $g$, we have
\be 
\int_{b}^{\infty} r^{d-1}g(r) dr \leq C_2 b^{d-1} e^{- c b}.
\lab{pe4aa}
\ee
Using this inequality in (\ref{pe4}) we get for all $k$ sufficiently large
\be
P(E_k)  \leq  C_3 k^a \log k \lb \hr_{n_k}(\ga)^{-1} \tL_{n_{k+1}}(\ep) \rb^{d-1} e^{- c \hr_{n_k}(\ga)^{-1} \tL_{n_{k+1}}(\ep)}.
\lab{pe4a}
\ee
Since $\f{n_k+1}{n_k} \to 1$ as $k \to \infty$, we have for all $k$ sufficiently large
\[ c \hr_{n_k}(\ga)^{-1} \tL_{n_{k+1}}(\ep) \geq (1+\ep)(1 - \del) \log n_k. \]
Substituting in (\ref{pe4a}), we get for all $k$ sufficiently large
\be
P(E_k)  \leq  C_4 k^a (\log k )^d \f{1}{k^{a(1+\ep)(1-\del)}} = C_4 \f{(\log k )^d}{k^{a((1+\ep)(1-\del) - 1)}}, 
\lab{pe4b}
\ee
which is summable since $a((1+\ep)(1-\del) - 1) > 1$. By the Borel-Cantelli lemma, almost surely, only finitely many of the $E_k$ and hence the $D_n(\hr_n(\ga), \tL_n(\ep))$ occur. In other words, almost surely, we have
\[ \f{L_n(\ga)}{\hr_n(\ga) \log n} \leq \f{1+\ep}{c}, \]
for all $n$ sufficiently large. (\ref{pe3}) now follows since $\ep > 0$ is arbitrary. 

To prove (\ref{pe3a}), we proceed as in the proof of (\ref{pe3}) above with $\tL_n(\ep)$ replaced by $\tL_n(b) = n^{-\f{b}{d}}$ where $b < \f{c-2d}{c-d}$ is arbitrary. It can be easily verified that in this case the bound for the integral on the left hand side in (\ref{pe4aa}) will be $C_2 b^{d-c}$. Using this bound in (\ref{pe4}) we get for all $k$ sufficiently large
\bea 
P(E_k) & \leq & C_3 n_{k+1}^2 \hr_{n_k}(\ga)^d \lb \hr_{n_k}^{-1} \tL_{n_{k+1}}(b) \rb^{d-c} \no \\
  & \leq & C_4 k^a \log k \lb \f{k^a}{\log k} \rb ^{1 - \f{c}{d}} k^{-ab \lb 1 - \f{c}{d} \rb } \no \\
 & \leq & C_4 \f{\lb \log k \rb ^{\f{c}{d}}}{k^{a\lb \f{c}{d} - 2 - b \lb \f{c}{d} - 1 \rb \rb}}. 
\lab{pe4c}
\eea
Since $b < \f{c-2d}{c-d}$, we have $\lb \f{c}{d} - 2 - b \lb \f{c}{d} - 1\rb \rb > 0$. Hence we can and do choose $a$ large enough so that $\sum P(E_k) < \infty$. By the Borel-Cantelli lemma, almost surely $L_n(\ga)^d \leq n^{-b}$ for all $n$ sufficiently large. In other words, 
\[ \limsup_{n \to \infty} \f{\log L_n(\ga)}{\log n} \leq - \f{b}{d}. \]
(\ref{pe3a}) now follows since $b < \f{c-2d}{c-d}$ is arbitrary. \qed

\begin{remark} While the above technique gives us tighter bounds, which are of independent interest, we will need a somewhat stronger statement on the behaviour of $L_n$ for the localization arguments we present below. Suppose we take $a = 1$ in the above proof, that is, we do not use the subsequence argument. Then following the above proof it is easy to see that
\be  \sum_{n=1}^{\infty} P\lb D_n(\hr_n(\ga), \tL_n(\ep) \rb < \infty, 
\lab{pe5}
\ee
holds in the case when $g(r) = e^{-rc}$ as $r \to \infty$ for some $c > 0$ if we take
\[ \tL_n(\ep) = c^{-1} (2 + \ep) \hr_n(\ga) \log n, \]
and in the case $g(r) = o(r^{-c})$ as $r \to \infty$ for some $c > 3d$, if we take $\tL_n(b) = n^{-\f{b}{d}}$ with $b < \f{c - 3d}{c-d}$. For instance, for (\ref{pe4c})  to be summable with $a=1$ we must have $\f{c}{d} - 2 - b\lb \f{c}{d} - 1 \rb > 1$, which holds if $b < \f{c-3d}{c-d}$.
\lab{rem3}
\end{remark}

\begin{prop} Let $d_n$ be as defined in (\ref{e1}). 
Suppose that $g(r) = o(e^{-cr})$ as $r \to \infty$ for some $c > d$. Then almost surely,
\be \liminf_{n \to \infty} \f{\al n d_n^d}{\log n} \geq 1. 
\lab{pe5a} 
\ee
If $g(r) = o(r^{-c})$ as $r \to \infty$ for some $c > 3d$, then almost surely 
\be \liminf_{n \to \infty} \f{\al n d_n^d}{\log n} \geq \f{c - 3d}{c - d}.
\lab{pe5b}
\ee
\lab{p6}
\end{prop}

{\bf Proof. }  We first prove (\ref{pe5a}). Fix $\ep \in (0,1)$ and choose $\del, \ep_1$ such that $0 < \del < \ep < \ep_1 < 1$, and let $\eta = 1- \ep$. Let $\tr_n^d = \f{\del \log n}{n \th}$, where $\th$ is the volume of the unit ball in $\mR^d$. Let 
\be \tL_n(\ep) = (2+\ep)c^{-1} \hr_n(\eta) \log n = \f{(2+\ep)}{c \al} \f{ \lb \log n \rb^{1 + \f{1}{d}}}{n^{\f{1}{d}}}, \qquad n \geq 1, 
\lab{pe6}
\ee
and define the sequence of events
\be D_n(\hr_n(\eta), \tL_n(\ep)) = \{ L_n(\eta) > \tL_n(\ep) \},
\lab{pe7}
\ee
where $L_n(\eta)$ is the longest edge in the graph $G_n(\hr_n(\eta))$. From Remark~\ref{rem3}, we have
\be \sum_{n=1}^{\infty} P \lb D_n(\hr_n(\eta), \tL_n(\ep)) \rb < \infty.
\lab{pe8}
\ee
Let $\ka_n$ be the packing number of $S$ by balls of radius $\tL_n(\ep_1)$. Then for sufficiently large $n$, we have
\be \ka_n \geq \f{C_1}{\tL_n(\ep_1)^d} = \f{C_2 \, n}{(\log n)^{d+1}}.
\lab{pe9}
\ee
Let $\{x_1^{(n)}, x_2^{(n)}, \ldots , x_{\ka_n}^{(n)}\}$ be a deterministic set of points in $S$ such that the balls $B(x_i^{(n)}, \tL_n(\ep_1))$, $i=1,2, \ldots ,\ka_n$ are disjoint. Let $E_n^{(i)}$ $\lb \tE_n^{(i)}\rb$ be the event that there is exactly one point of $\Pn$ in the ball $B(x_i^{(n)}, \tr_n)$ which is isolated (respectively, has no edge to any point of $\Pn \cap B(x_i^{(n)}, \tL_n(\ep_1))$) in the graph $G_n(\hr_n(\eta))$.

Suppose we show that 
\be \sum_{n=1}^{\infty} P \lb \cap_{i=1}^{\ka_n} \lb E_n^{(i)} \rb^c \rb < \infty.
\lab{pe10}
\ee
It then follows by the Borel-Cantelli lemma that almost surely, for all $n$ sufficiently large, $E_n^{(i)}$ occurs for some $i$. In other words, there is a vertex in $G_n(\hr_n(\eta))$ which is isolated, that is, $d_n > \tr_n(\eta)$. The result will then follow since $\eta < 1$ is arbitrary. Thus it remains to show (\ref{pe10}).

\be P \lb \cap_{i=1}^{\ka_n} \lb E_n^{(i)} \rb^c \rb \leq P \lb \cap_{i=1}^{\ka_n} \lb E_n^{(i)} \rb^c \cap D_n^c \rb +
P(D_n),
\lab{pe11}
\ee
where $D_n = D_n(\hr_n(b), \tL_n(\ep))$ is as defined in (\ref{pe7}). (\ref{pe10}) will follow from (\ref{pe8}) and (\ref{pe11}) provided we show that
\[ \sum_{n=1}^{\infty} P \lb \cap_{i=1}^{\ka_n} \lb E_n^{(i)} \rb^c \cap D_n^c \rb < \infty. \]
From (\ref{pe6}) we have 
\be \tL_n(\ep) + \tr_n = \f{ \lb \log n \rb^{1 + \f{1}{d}}}{n^{\f{1}{d}}} \lb \f{(2+\ep)}{c \al} + \f{\del^{\f{1}{d}}}{\log n} \rb \leq \tL_n(\ep_1),
\lab{pe11a}
\ee
for all $n$ sufficiently large. It follows that
\bea
P \lb \cap_{i=1}^{\ka_n} \lb E_n^{(i)} \rb^c \cap D_n^c \rb & \leq & P \lb \cap_{i=1}^{\ka_n} \lb \tE_n^{(i)} \rb^c \rb = \prod_{i=1}^{\ka_n} P \lb \lb \tE_n^{(i)} \rb^c \rb \nn \\
 & \leq & \exp \lb - \ka_n P \lb \tE_n^{(1)} \rb \rb,
\lab{pe12}
\eea
where we have used the fact that the events $\tE_n^{(i)}$, $i=1,2, \ldots ,n,$ are independent and the inequality $1-x \leq e^{-x}$ .  Let $A_n = B(x_n^{(1)}, \tL_n(\ep_1)) \setminus B(x_n^{(1)}, \tr_n)$. Then using Lemma~\ref{l1} we get
\bea
P \lb \tE_n^{(1)} \rb & = & n \th \tr_n^d \exp \lb - n \th \tr_n^d \rb \f{1}{n \th \tr_n^d} \int_{B(x_n^{(1)}, \tr_n)}
\exp \lb -n \int_{A_n} g \lb  \f{|z-y|}{\hr_n(\eta)} \rb dz \rb dy \nn \\
 & \geq & C_3 \f{\log n}{n^{\del}} \exp \lb - n \hr_n(\eta)^d \al \rb \nn \\
 & = & C_3 \f{\log n}{n^{\del + \eta}},
\lab{pe13}
\eea
where the inequality in the second line above is obtained by first making the standard change of variable and then replacing $\hr_n^{-1}A_n$ by $\mR^d$.
Substituting from (\ref{pe9}) and (\ref{pe13})  in (\ref{pe12}) we get
\beas
P \lb \cap_{i=1}^{\ka_n} \lb E_n^{(i)} \rb^c \cap D_n^c \rb & \leq & \exp \lb - C_4 \f{n}{(\log n)^{d+1}} \f{\log n}{n^{\del+\eta}} \rb \nn \\
 & = & \exp \lb - C_4 \f{n^{\ep - \del}}{(\log n)^{\del}} \rb,
\eeas
which is summable since $\ep > \del$. This completes the proof of (\ref{pe5a}).

To prove (\ref{pe5b}), we proceed exactly as in the proof of (\ref{pe5a}) with some minor changes.  Choose $\eta, b, b_1$ such that $0 < \eta < b_1 < b < \f{c-3d}{c-d}$. Replace $\tL_n(\ep)$ by $\tL_n(b) = n^{-b/d}$ and $\tL_n(\ep_1)$ by $\tL_n(b_1) = n^{-b_1/d}$.  Choose $\del > 0$ so that $\del + \eta < b_1.$ The inequality (\ref{pe9}) will change to
\be \ka_n \geq C_1 n^{b_1}.
\lab{pe13a}
\ee
The analog of (\ref{pe11a}) holds since $b_1 < b < 1$. The bound obtained in (\ref{pe13}) holds in this case as well. Thus, substituting from (\ref{pe13}), (\ref{pe13a}) we get
\beas 
P \lb \cap_{i=1}^{\ka_n} \lb E_n^{(i)} \rb^c \cap D_n^c \rb & \leq & \exp \lb - C_4 n^{b_1} \f{\log n}{n^{\del + \eta}} \rb \nn \\
 & = & \exp \lb - C_4 n^{b_1 - \del - \eta} \log n  \rb,
\eeas
which is summable since $b_1 > \del + \eta$. This proves (\ref{pe5b}). \qed

{\bf Proof of Theorem~\ref{t1}. } Theorem~\ref{t1} now follows from Propositions~\ref{p5} and \ref{p6}. \qed

{\bf Proof of Theorem \ref{t2}. } Let $\beta$ be as defined in (\ref{e3b}). Fix $\ga > \beta$ and choose $b > 1$ such that $\ga g \lb \f{bc}{\ga \th} \rb > 1.$ Recall that $\hr_n^d(\ga) = \f{\ga \log n}{n \al}$ and $\tr_n(b) = \f{b\log n}{n \th}$. 
Consider the graph $G_n(\hr_n(\ga))$. A one hop path is said to exist between $X_i,X_j \in \Pn$ if there is a vertex $X_k \in \Pn$ such that $\{X_i,X_k\}$ and $\{X_k,X_j\}$ are edges in $G_n(\hr_n(\ga))$. Let $E_n = E_n(\ga)$ be the event that there is a vertex $X \in \Pn$ such that $X$ does not have a one hop path to some point in $\Pn \cap B(X, \tr_n(b))$ in the graph $G_n(\hr_n(\ga))$. 

Suppose we show that $P(E_n) \to 0$ as $ n \to \infty$. It follows that every vertex in $G_n(\hr_n(\ga))$ is connected to each one of its $\tr_n(b)$ neighbour via a one hop path {\it whp}. However existence of a path to all the $\tr_n(b)$ neighbour implies that the graph is connected {\it whp} (Theorem 13.7, \cite{Penrose03}). This will prove the result. To this end we estimate $P(E_n)$.

For any $X \in \Pn$, let $E_n(X)$ be the event that $X$ is not connected to at least one vertex in $\Pn \cap B(X, \tr_n(b))$ in the graph $G_n(\hr_n(\ga))$ via a one hop path. Then by the Campbell-Mecke formula we get
\[ P(E_n) \leq E \left[ \sum_{X \in \Pn} 1_{E_n(X)} \right] = n P^o \lb E_n(O) \rb .\]
Let $H_n$ be the event that the number of points of $\Pn$ in $B(O,\tr_n(b))$ is does not exceed $a \log n$. Then by  Theorem 6.14, \cite{Penrose03}, we can and do choose $a$ sufficiently large so that $P(H_n^c) \to 0$ as $n \to \infty$ . Let $A_n$ be the event that the origin is not connected to a point chosen uniformly at random in the ball $B(O,\tr_n(b))$ via a one hop path in the graph $G_n(\hr_n(\ga))$. Since $P^o(E_n(O)) \leq P^o(E_n(O) \cap H_n) + P^o(H_n^c)$, by the above observation and using the union bound, it suffices to show that $n (\log n) P^o(A_n) \to 0$ as $n \to \infty$.
\be
P^o(A_n)  =  \f{1}{\th \tr_n(b)^d} \int_{B(O, \tr_n(b))} \exp \lb - n \int_S g \lb \f{d(O,y)}{\hr_n(\ga)} \rb
g \lb \f{d(y,z)}{\hr_n(\ga)} \rb dz \rb dy.
\lab{pe14}
\ee
The integrand in (\ref{pe14}) is the probability that there is no point in $\Pn \setminus \{O,y\}$ that connects to both the origin and $y$ in $G_n(\tr_n(\ga))$ under $P^{o,y}$. Changing the variables from $(y,z)$ to $(u,v)$ with
\[ u = \f{y}{\hr_n(\ga)} \qquad \mbox{ and } \qquad  v = \f{z}{\hr_n(\ga)}, \]
in (\ref{pe14}), and using Lemma~\ref{l1} we get
\be 
P^o(A_n)  \leq C_2  \f{\hr_n(\ga)^d}{\th \tr_n(b)^d} \int_{B(O, \hr_n(\ga)^{-1} \tr_n(b))} \exp \lb - \al n \hr_n(\ga)^d  g(|u|) \rb du.
\lab{pe15}
\ee
Since $\hr_n(\ga)^{-1} \tr_n(b) = \f{b \al}{\ga \th},$ and $g$ is non-increasing, (\ref{pe15}) can be bounded by
\[ P^o(A_n)  \leq C_2  \exp \lb - \ga g \lb \f{b \al}{\ga \th} \rb \log n \rb.
\]
Consequently
\[ n (\log n) P^o(A_n)  \leq C_2 \f{\log n}{n^{\ga g \lb \f{b \al}{\ga \th} \rb - 1}} \to 0, \]
since $\ga g \lb \f{b \al}{\ga \th} \rb > 1$. \qed

{\bf Proof of Proposition \ref{p2}. } Recall that $\hr_n = \f{\ga \log n}{\al n}$ and $\br_n = a_n \hr_n$, where the sequence $\{a_n\}$ satisfies $G(a_n) \to e^{-\beta}$ as $n \to \infty$ with $G$ as defined in (\ref{e4a}). The points of $\Pn$ under $P^o$ located in $S \setminus B(O,\br_n)$ that have an edge to the origin form a non-homogenous Poisson point process and hence
\beas
P^o \lb L_n^o \leq \br_n \rb & = & \exp \lb - n \int_{S \setminus B(O,\br_n)} g \lb \f{|z|}{\hr_n} \rb dz \rb \\
 & = & \exp \lb - n \hr_n^d \int_{\hr_n^{-1} \lb S \setminus B(O,\br_n) \rb} g(|z|) dz \rb \\
 & = & \exp \lb - n \hr_n^d \int_{ \mR^d \setminus B(O, \hr_n^{-1} \br_n) } g(|z|) dz \rb 
\exp \lb - n \hr_n^d \int_{\mR^d \setminus \hr_n^{-1} S } g(|z|) dz \rb \\
 & \to & e^{-e^{-\beta}},
\eeas
as $n \to \infty$ by Lemma~\ref{l1}. \qed

Since Proposition~\ref{p3} has already been proved, we now prove Proposition~\ref{p4}.

{\bf Proof of Proposition \ref{p4}. } By the Campbell-Mecke formula,
\beas
E \left[ \f{D_n(k_n)}{n} \right] & = &  P^o \lb \mbox{deg(O) } \geq k_n \mbox{ in } G_n(\hr_n) \rb \\
 & = & P\lb Z_n \geq k_n \rb, 
\eeas
where $Z_n$ is a Poisson random variable with mean
\beas 
m_n & = & n \int_S g \lb \f{|z|}{\hr_n} \rb dz \\
 & = & n \hr_n^d \int_{\hr_n^{-1}S} g(|z|) dz \sim \al n \hr_n^d.
\eeas
The result now follows from the central limit theorem. \qed

{\bf Proof of Theorem \ref{t3}. } We first prove (\ref{e8}). Fix $\ep > 0$ and choose $\del \in (0, \ep)$ such that $\eta = (1-\del)(1+\ep) - 1 > 0$. Let $n_k = k^a$ where $a > 1$ is chosen such that $a \eta > 1$. Define the sequence
\[ c_n = (1+\ep) \ga H_+^{-1}\lb \f{1+\ep}{\ga} \rb \, \log n . \]
Then using the coupling and the union bound we get
\bea
P \lb \cup_{n=n_k}^{n_{k+1}} \lb \Del_n(\hr_n) \geq c_n \rb \rb & \leq & P \lb \Del_{n_{k+1}}(\hr_{n_k}) \geq c_{n_k} \rb \nn \\
& \leq & n_{k+1} P^o \lb \mbox{deg}(O) \geq c_{n_k} \mbox{ in } G_{n_{k+1}}(\hr_{n_k}) \rb.
\lab{pe16}
\eea
The degree at the origin under $P^o$ in $G_{n_{k+1}}(\hr_{n_k}(\ga))$ is Poisson distributed with mean 
\be
m_k = n_{k+1} \int_S g \lb \f{|z|}{\hr_{n_k}(\ga)} \rb dz,
\lab{pe17}
\ee
and hence using the Chernoff bound (see Lemma 1.2, \cite{Penrose03}) we have
\be 
P^o \lb \mbox{deg}(O) \geq c_{n_k} \mbox{ in } G_{n_{k+1}}(\hr_{n_k}) \rb
\leq \exp \lb - m_k H \lb \f{c_{n_k}}{m_k} \rb \rb.
\lab{pe18}
\ee
By Lemma \ref{l1}, $m_k \sim n_{k+1} \hr_{n_k}^d(\ga) \al = \ga \f{n_{k+1}}{n_k} \log n_k$. Since $\f{n_{k+1}}{n_k} \to 1$ as $k \to \infty$ we have for all $k$ sufficiently large
\be
 (1 - \del) \ga \log n_k \leq m_k \leq (1+\del) \ga \log n_k.
\lab{pe19}
\ee
Since $\del < \ep$, we have for all $k$ sufficiently large,
\be 
\f{c_{n_k}}{m_k} \geq \f{1+\ep}{1+\del} H_+^{-1}\lb \f{1+\ep}{\ga}\rb \geq H_+^{-1}\lb \f{1+\ep}{\ga} \rb.
\lab{pe20}
\ee
Substituting from (\ref{pe19}), (\ref{pe20}) in (\ref{pe18}) and using the fact that $H$ is increasing in $[1, \infty)$ we get
\bea
P^o \lb \mbox{deg}(O) \geq c_{n_k} \mbox{ in } G_{n_{k+1}}(\hr_{n_k}) \rb
 & \leq & \exp \lb - (1-\del) \ga \log n_k \lb \f{1+\ep}{\ga} \rb \rb \nn \\
 & \leq & \exp \lb - (1 + \eta) \log n_k \rb. 
\lab{pe21}
\eea
Substituting from (\ref{pe21}) in (\ref{pe16}) we get
\bea
P \lb \cup_{n=n_k}^{n_{k+1}} \lb \Del_n \geq c_n \rb \rb & \leq &
n_{k+1} \exp \lb - (1 + \eta) \log n_k \rb \nn \\
 & \leq & C_1 \f{1}{k^{a \eta}},
\lab{pe21a}
\eea
which is summable in $k$. Hence by the Borel-Cantelli lemma, almost surely, only finitely many of the events
$\cup_{n=n_k}^{n_{k+1}} \lb \Del_n \geq c_n \rb $, and hence, only finitely many of the events
$\{ \Del_n \geq c_n \}$ occur. It follows that almost surely,
\[ \f{\Del_n}{\log n} \leq (1 + \ep) \ga H_+^{-1}\lb \f{1+\ep}{\ga} \rb, \]
for all $n$ sufficiently large. Since $\ep > 0$ is arbitrary and $H_+^{-1}$ is increasing, we have that almost surely
\[ \limsup_{n \to \infty} \f{\Del_n}{\log n} \leq \ga H_+^{-1}\lb \f{1}{\ga} \rb. \]
This proves (\ref{e8}). To prove (\ref{e9}) it suffices by (\ref{e8}) to show that 
\be 
\liminf_{n \to \infty} \f{\Del_n}{\log n} \geq \ga H_+^{-1}\lb \f{1}{\ga} \rb.
\lab{pe21b}
\ee
The proof of (\ref{pe21b}) is similar to that of (\ref{pe5a}) and we will borrow notations used in proving (\ref{pe5a}). Fix $\ep_1 \in (0,1)$. Since $H_+^{-1}$ is increasing we can and do choose $\ep \in (0, \ep_1)$ and $\eta > 0$ such that
\be
(1-\eta)^{-1} H_+^{-1} \lb \f{1 - \ep_1}{\ga} \rb \leq H_+^{-1} \lb \f{1 - \ep}{\ga} \rb.
\lab{pe21c}
\ee
Let $\tL_n(\ep)$ be as defined in (\ref{pe6}) and let $D_n$ defined as in (\ref{pe7}). Let $\ka_n$ be the packing number of $S$ by balls of radius $\tL_n(\ep_1)$ which satisfies the inequality in  (\ref{pe9}) for all $n$ sufficiently large.  Let $\{x_1^{(n)}, x_2^{(n)}, \ldots , x_{\ka_n}^{(n)}\}$ be a deterministic set of points in $S$ such that the balls $B(x_i^{(n)}, \tL_n(\ep_1))$, $i=1,2, \ldots ,\ka_n,$ are disjoint. Define the sequence $\{c_n\}_{n \geq 1}$ by
\be
c_n =  \ga H_{+}^{-1} \lb \f{1-\ep_1}{\ga} \rb \log n. 
\lab{pe22}
\ee
Recall that $\tr_n = \f{\del \log n}{n \th}$, where $\th$ is the volume of the unit ball in $\mR^d$. Choose $\del < \ep$ sufficiently small so that 
\be 
m_n(y) = n \int_{A_n} g \lb  \f{|z-y|}{\hr_n(\ga)} \rb dz,\lab{pe22a}
\ee
where $A_n = B(x_n^{(1)}, \tL_n(\ep_1)) \setminus B(x_n^{(1)}, \tr_n)$ satisfies
\be
\f{m_n}{n \hr_n(\ga)^d \al} \geq 1 - \eta,
\lab{pe22b}
\ee
for all $n$ sufficiently large. Such a choice of $\del$ is possible since
\[ \f{m_n}{n \hr_n(\ga)^d} \to \int_{\mR^d \setminus B \lb O, \lb \f{\del \al}{\ga \th} \rb^{\f{1}{d}} \rb } g(|z|) dz,\]
as $n \to \infty$ and the right hand side expression in the above equation converges to $\al$ as $\del \to 0$. 
 
Let $E_n^{(i)}$ $\lb \tE_n^{(i)}\rb$ be the event that there is exactly one point of $\Pn$ in the ball $B(x_i^{(n)}, \tr_n)$ whose degree is at least $c_n$ (which has an edge to at least $c_n$ many points in $\Pn \cap B(x_i^{(n)}, \tL_n(\ep_1))$) in the graph $G_n(\hr_n(\ga))$. (\ref{pe21b}) follows if we show (\ref{pe10}) by an argument similar to the one that appears below it. Following the same steps as in the proof of (\ref{pe10}), it suffices to show that $\exp \lb - \ka_n P \lb \tE_n^{(1)} \rb \rb$ is summable. 
\be
P \lb \tE_n^{(1)} \rb = n \th \tr_n^d \exp \lb - n \th \tr_n^d \rb \f{1}{n \th \tr_n^d} \int_{B(x_n^{(1)}, \tr_n)} \exp \lb - m_n(y) H \lb \f{c_n}{m_n(y)} \rb \rb dy,
\lab{pe23}
\ee
where $m_n$ is as defined in (\ref{pe22b}). By the choice of $\del$ and (\ref{pe21c}) we have for all $n$ sufficiently large
\[ \f{c_n}{m_n} \leq (1 - \eta)^{-1} H_{+}^{-1} \lb \f{1-\ep_1}{\ga} \rb \leq H_+^{-1} \lb \f{1 - \ep}{\ga} \rb. \]
By the standard change of variable, it is easy to see that $m_n \leq n \hr_n^d(\eta) \al = \ga \log n$. Using these two inequalities in (\ref{pe23}) and the fact that $H(x)$ is increasing for $x > 1$ we have for all $n$ sufficiently large
\beas 
P \lb \tE_n^{(1)} \rb & \geq & C_3 \f{\log n}{n^{\del}} \exp\lb - (1 - \ep) \log n \rb \\
 & = & C_3 \f{\log n}{n^{1 + \del - \ep}}.
\eeas
It follows from (\ref{pe9}) and the above inequality that
\bea 
\exp \lb - \ka_n P \lb \tE_n^{(1)} \rb \rb & \leq & \exp \lb - C_4 \f{n}{(\log n)^{d+1}} \f{\log n}{n^{1 + \del - \ep}} \rb \nn \\
 & = & \exp \lb - C_4 \f{n^{\ep - \del}}{(\log n)^{d}} \rb,
\lab{pe24}
\eea
which is summable since $\ep > \del$. This proves (\ref{e9}) since $\ep_1 > 0$ is arbitrary.

To prove (\ref{e10}) we proceed as in the proof of (\ref{e9}) with the modifications similar to those in the proof of  (\ref{pe5b}). Fix $b_1 \in \lb 0, \f{c-3d}{c-d} \rb$ and choose $b \in \lb b_1 , \f{c-3d}{c-d} \rb$. Choose $\ep, \eta > 0$ so that 
\[ (1-\eta)^{-1} H_+^{-1} \lb \f{b_1}{\ga} \rb \leq H_+^{-1} \lb \f{b_1 - \ep}{\ga} \rb. \]
Take $c_n = \ga H_+^{-1} \lb \f{b_1}{\ga} \rb \log n$. Replace $\tL_n(\ep), \tL_n(\ep_1)$ by $\tL_n(b)$ and $\tL_n(b_1)$ respectively, in the proof of (\ref{e9}) above. Choose $\del < \ep$ so as to satisfy (\ref{pe22b}). Since $\ka_n$ now satisfies (\ref{pe13a}), proceeding as in the proof of (\ref{e9}) we will end up with the same bound as in (\ref{pe24}) which is summable since $\ep > \del$. The result now follows since $b_1 < \f{c-3d}{c-d}$ is arbitrary. \qed

{\bf Proof of Theorem \ref{t4}. } Part $(i)$ of the proof follows immediately from Theorem~\ref{t1}~$(ii)$. The rest of the proof is entirely analogous to the proof of Theorem~\ref{t3} with the obvious changes. We illustrate this by proving part (ii), the proof of which is similar to that of (\ref{e8}).

Fix $\ep > 0$ so that $\f{1+\ep}{\ga} < 1$. Let $\del, \eta, a, n_k$,  be as in the proof of (\ref{e8}). Define $c_n$ as
\[ c_n = (1-\ep) \ga H_-^{-1} \lb \f{1+\ep}{\ga} \rb \,  \log n . \]
Consider the sequence $s_k = P\lb \cup_{n=n_k}^{n_{k+1}} \lb \del_n (\hr_n(\ga)) \leq c_n \rb \rb \leq P \lb \del_{n_k}(\hr_{n_{k+1}}) \leq c_{n_{k+1}} \rb$. Proceeding as in the proof of (\ref{e8}) and using the fact that $H$ is decreasing in $[0,1]$ we obtain the same bound for $s_k$ as on the right hand side of (\ref{pe21a}). It then follows that the sequence $s_k$ is summable from which (\ref{e11}) follows. \qed

\end{document}